\documentclass[10pt,a4paper]{article}
\usepackage[utf8]{inputenc}
\usepackage{amsmath}
\usepackage{amsfonts}
\usepackage{amssymb}
\usepackage{graphicx}
\usepackage{authblk}
\usepackage{hyperref}
\usepackage{color}
\usepackage{mathtools}
\usepackage[ruled,vlined]{algorithm2e}
\usepackage{algorithmic}
\SetKwInput{KwInput}{Input}                
\SetKwInput{KwOutput}{Output}              

\author[1]{Rajeev Kumar\\ Society for Natural Technology Research,\\ 
Dept. of IT\&E, Governemnt of West Bengal, India.\\  
rajeev.ips@nltr }
\date{}

\providecommand{\keywords}[1]
{
  \small	
  \textbf{\textit{Keywords --- }} #1
}
\begin{document}
\title{   Finding discrete logarithm in $F_p^* $ }
\maketitle
\newcommand{\cmnt}[2]{}

\begin{abstract}
Difficulty of calculation of discrete logarithm for any arbitrary Field is the basis for security of several popular cryptographic solutions.  Pohlig-Hellman method is a popular choice to calculate discrete logarithm in finite field $F_p^*$. Pohlig-Hellman method does yield good results if p is smooth ( i.e. p-1 has small prime factors).  We propose a practical  alternative to Pohlig-Hellman algorithm for finding discrete logarithm modulo prime. Although, proposed method, similar to Pohlig-Hellman reduces the problem to group of orders $p_i$ for each prime factor and hence in worst case scenario (when p=2q+1 , q being another prime) order of run time remains the same. However in proposed method, as  there is no requirement of combining the result using Chinese Remainder Theorem and do other associated work ,run times are much faster.   \\
\end{abstract}

\section{Introduction}
Discrete Logarithm problem in field $F_p^*$ , where p is prime ( $ g,\beta \in F_p^*$ ), where $g$ a generator, can be defined as finding x for which $g^x \equiv \beta \pmod {p} $. \\ Traditional methods avaialble for finding discrete logarithm  are `baby steps giant steps' \cite{babystepgiantstep} , Silver Pohlig-Hellman \cite {silver-pohilg} and Pollard's rho/Kangaroo Method \cite {pollard} \cite{rho}. Silver-pohlig-hellman being the fastest if prime is smooth, i.e prime factors of p-1 can be easily found and such factors are far smaller then p. Run time of pohlig-hellman method in worst case is  $\mathcal{O}(\sqrt{p})$, however what is interesting is that it is a lot more efficient when the order of p is smooth. A smooth order means that the order has prime factors, which are far smaller. If the prime factorisation of $p-1 = p_1^{e_1}p_2^{e_2}\dot{...} p_k^{e_k}$ the complexity is $ \mathcal{O} {(\sum_i^k e^i(\log{p}+ \sqrt{p_i}))} $ 
Our proposed method uses the idea of calculating the $p_i$th root  for all  prime factors which eventually leads us to the desired result. The result is obtained directly without the need of combining results for all $p_i^{e_i} $ through Chinese remainder theorem.
The run time complexity of proposed method ,  for prime p whose order is B-smooth is $ \mathcal{O} {(\sum_i^k e^i( \sqrt{ p_i}))}$ , almost same as that of Pohlig-Hellman method. However in practical terms absence of extra work done in Pohlig-Hellman method makes this a better candidate for solving discrete logarithm problem , although there is no theoretical change in the upper bound of the run-time complexity.
\section{Basic Background}
Let $g,\beta \in F^*_p $ where $g$ is  a generator and $ p-1= p_1^{e_1}p_2^{e_2}\dot{...} p_k^{e_k}$. Now by definition , of $g$ being generator , we can always find $x \in \{0,..,p-1\}$ such that $ g^x \equiv \beta$.\\
\\Also , by Fermat's Theorem $ g^{(p-1)} \equiv 1$. 
In $F^*_p$ , if $i^{th}$ root of any element $\beta$ exists, there are exactly $i$ such roots. Also all such i roots $ \{ \beta_1,\beta_2,..\beta_i\}$ form a cyclic group. Other elements of this group can be generated from $\beta_1$ by successively multiplying it with  $g ^{(p-1)/i}$ for any $g$, a generator $ \in F^*_p $. However 

\section{Algorithm }
Given in field $F_p^*$ , where p is prime ( $ g,\beta \in F_p^*$ ), where $g$ a generator and $ p-1= p_1^{e_1}p_2^{e_2}\dot{...} p_k^{e_k}$.
Let $\sum_i^k e_i = s(say)$.
First we note that in s  steps of exponentiation  we get $\{\beta,\beta^{p_1}, \dot{...},\beta^{ p_1^{e_1}p_2^{e_2}\dot{...} p_k^{e_{k-1}}}, 1\}$ or $\{\beta, \beta_1,\beta_2,..\beta_{s-1},1\}$. We note that $\beta^{p_1^{e_1}p_2^{e_2}\dot{...} p_k^{e_k}} \equiv 1$. We can write these $s+1$ elements as trace back list.$\{1, \beta_{s-1}, \beta_{s-2},\beta_{s-3},..,\beta_1,\beta \}$ . We know by Fermat's Theorem that $g^{(p-1)}\equiv 1 $. Now idea is to keep taking $p_i^{th}$ , root starting from $g^{(p-1)}$ , so that we trace back the list generated from $\beta$. Basic idea is to generate a list equivalent to trace back list where elements at each positions are equal , by taking the appropriate $p_i^{th} $ root .\\
i.e.  $ \{g^{p-1}, g^{x_1}, g^{x_2},..,g^x \} \equiv   \{1, \beta_{s-1}, \beta_{s-2},\beta_{s-3},..,\beta_1,\beta \}  $, trace back list.\\ $x_1,x_2,..x $ are obtained by taking successive $ p_i^{th}$ root starting from $g^{p-1}$\\
We note that one of the  $p_i^{th}$ root of $ g^{(p-1)}$ would be given $g^{(p-1)/p_i}$ and it can be clearly seen that such root exist as $p_i\text{ } |\text{ } (p-1)$.Now we start with by taking $p_k$ root of $ g^{p-1} \equiv 1$ .Let $g^{(p-1)/p_k} \equiv \delta $. Now  $ \beta_{s-1} $ would be member of the group of order $p_k$ i.e. \\ $ \beta_{s-1} \in \{ \delta, \delta g^{(p-1)/p_k},\delta g^{2(p-1)/p_k},..\delta g^{(p_k-1)(p-1)/p_k} \} $  \\ or \\  $ \beta_{s-1} \in \{\delta, \delta \epsilon_0,\delta \epsilon_0^2,..\delta \epsilon_0^{p_k-1} \} $  where $ \epsilon_0 \equiv g^{(p-1)/p_k} $. \\
\\$\beta_{s-1} $ can be located in the cyclic group of order $p_k$ through any of the baby step - giant steps method or pollard's rho method in $\mathcal{O}{(\sqrt{p_k})}$ time.\\
Following this first step we will get $g^r \equiv \beta_{s-1}$ where $ r$ is multiple of $p_1^{e_1},p_2^{e_2},..,p_k^{e_k-1}$. We can continue taking $p_i^{th}$ root of $g^r$ and we retrace the path from $\beta_{s-1}$, $\beta_{s-2} ,..\beta_1$  till we reach $ g^x \equiv \beta$.

\subsection{Pseudo Code}
Please note pseudo code is for illustrative purposes only.

\begin{algorithm}[H]
\caption{Pseudo Code  to calculate x such that $ g^x \equiv \gamma $ i.e. $\log{\gamma}$}
\KwInput{$g ,\gamma \in F_p^* , p$; $g $ is a generator and $ p = \prod_i^kp_i^{e_i}+1 $ where $p_i$ are all primes.}
\KwOutput{$x \text{ which is a solution of } g^x\equiv \gamma \pmod p $ }

\begin{algorithmic}
\STATE $Stack\_S \gets \gamma, elem \gets \gamma  $ 
\STATE $plist \gets [p_1,p_1,..p_k,p_k] \text{ i.e. list of all prime factors ,as many times as they appear.} $  
\WHILE{$ plist \text{ is not empty } $}
	\STATE $ elem \gets elem^{pop(plist)} , Stack \gets elem $ 
\ENDWHILE
\STATE $ gpow \gets p-1, plist \gets [p_1,p_1,..p_k,p_k] $
\WHILE{ $ plist \text{ is not empty }$}
	\STATE $ elem \gets pop(plist) , k \gets 0 $
	\STATE $ gpow \gets  gpow/elem , \delta \gets g^{gpow}$
	\IF {$\delta != \text{top of the Stack\_S } $ }
		\STATE Find $k$ such that $ \delta g^{k(p-1)/elem}\equiv \text{ Top of the stack}$ , through traditional methods 
	\ENDIF
	\STATE $ \text{Remove top element of Stack\_S} $
	\STATE $gpow \gets gpow + k*(gpow/elem)$
\ENDWHILE
\RETURN{$gpow$}
\end{algorithmic}
\end{algorithm}

\subsection{Example of algorithm's working}

We note $41=2^3.5^1+1 $ i.e list of prime factors of $p-1$ counting each appearance is  $ = [2,2,2,5] $ . Let us assume that $g=13$ is given as generator of the group and we have to find the discrete logarithm of 8, i.e we have to find $x$ such that $13^x \equiv 8 \pmod{41} $.\\
\\Now using the list of prime factors we obtain the trace back list $ [8, 8^2,8^{4},8^{8},8^{40}]$ or $[ 8,23,37,16,1]$ , all values calculated in field $ F^*_{41}$. So our trace back list is $[ 8,23,37,16,1]$. \\
\\ We start with generator now , we know $ 13^{40} \equiv 1$ \\
\\We take off the top element of trace back list , which is equal to $13^{40} $, now trace back list is $[ 8,23,37,16]$. 
\\Taking the fifth root , we get $ 13^{40/5} \equiv 13^{8} \equiv 10$
\\ All possible $5^{th}$ root of unity in field $ F^*_{41}$ roots can be written as $ \{13^8, 13^8*13^8,13^8*13^{16}, 13^8*13^{24},13^8*13^{32}\}$ or $[10, 18, 16, 37, 1] $. These are all possible $5^{th}$ root of unity in field $ F^*_{41}$ , hence top of the trace back list $16$ must be member of this group. we get $13^{24} \equiv 16 $. Locating this can be done by traditional methods like `baby steps giant steps' \cite{babystepgiantstep}  and Pollard's rho/Kangaroo Method \cite {pollard} \cite{rho} in $ \mathcal{O}(\sqrt{p_i})$, where $p_i$ is the prime factor and order of the group.
\\We take off the top element of trace back list , which is equal to $13^{24} $, now trace back list is $[ 8,23,37]$
\\ Backtracking,  we now take square root , we $ 13^{24/2} \equiv 13^{12} \equiv 4 $\\
All possible square roots are $ [ 13^{12} , 13^{12}*13^{20}]$ or $[4,37]$ . As our top of the trace back is $37$ we get $ 13^32 \equiv 37$.
\\
\\We take off the top element of trace back list ,which is equal to $13^{32} $, now trace back list is $[ 8,23]$
 Now again backtracking , we now take square root , we get $ 13^{32/2} \equiv 13^{16} \equiv 18 $\\.As our top of the trace back is $18$. All possible square roots are $ [ 13^{16} , 13^{16}*13^{20}]$ or $[18,23]$ . As element we are looking from our trace back list is 23 , we get $13^{36} \equiv 23 $ 
\\We take off the top element of trace back list ,which is equal to $13^{36} $, now trace back list is $[ 8]$
\\ Now for last time  backtracking , we again take square root , we get $ 13^{36/2} \equiv 13^{18} \equiv 8$. Our desired result  would definitely have been member of  $[ 13^{18} , 13^{18}*13^{20}] $ or $[8,33]$.  Here we get $ 13^{18} \equiv 8$ , so discrete $\log_{13}{8} \equiv 18  $ in $F_{41}^*. $\\

\section{Run-time Analysis}
We see trace back list, in field $F^*_p$  can be generated  $ \mathcal{O}( \log {p})$. taking $p_i^{th}$ root and locating the exact root we are looking as per our trace back list can be done in  $ \mathcal{O}(\sqrt {p_i})$ by baby-step-giant step algorithms for one particular $ p_i$. However this step has to be repeated for $e_i$ time for each $p_i$, hence total runtime for proposed method is $ \mathcal{O}(\log{p}+ \sum_i^k e_i\sqrt {p_i})$, which is bounded by same limits as provided by pohlig-hellman method.
However , in practical terms, no combination of various group results through CRT is to be done and other lesser book keeping during the algorithm leads to much faster running time.
\subsection{Run time comparison}
Python code implementation of both the proposed method and Pohlig-Hellman method has been provided for comparison. Test run were run 10000 times for each example and averaged time in micro seconds is given. It is easy to see that considerable practical advantages are offered by proposed method.
\begin{center}
 \begin{tabular}{||c c c c||} 
 \hline
       & Trial 1 & Trial 2& Trial 3 \\ [0.5ex] 
 \hline
 $prime\ p =$& 41 & 8101 & 200560490131 \\ 
 \hline
 $generator\ g=$ & 13 & 6 & 79 \\
 \hline
 $ h= $ & 8 & 7531 & 23 \\
 \hline
 $\log{h} \in F_p^* =$ & 18 & 6689 & 127013812855\\
 \hline\hline
 Dlp time & 16.3 & 47 & 274 \\ 
 \hline
 Pohlig-hellman  & 33.2 & 87.2 & 475 \\ [1ex] 
 \hline
\end{tabular}
\end{center}

\section{Advantages of Proposed Algorithm }
Before we count the advantages , let us explore why proposed method gives better run-times compared to pohlig-hellman.

\subsection{Why faster run-time ?}
Pohlig-Hellman method calculates discrete log for each prime factor $p_i^{e_i}$ and combines the result so obtained through Chinese remainder theorem.
However proposed method keeps taking $p_i^{th}$ root starting from $ g^{p-1}$ and every level it selects one out of $p_i$ values depending upon the trace back list already prepared by raising $\beta $, whose discrete logarithm is sought. This way $g^{p-1}$  reduces to $g^x$ mimicking  $\beta^{p-1}$'s journey to $\beta$. And $x$ is the desired result. Taking advantage of the property of $\beta$ a direct approach to calculating discrete Logarithm lends to faster run-time.

\subsection{Advantages}
1.  Advantages of proposed algorithm is more on practical side as it provides faster calculation of discrete logarithm put does not improve the worst case performance bound of Pohlig-Hellman method.
2. Pohlig-Hellman's run time  is agnostic to value whose log is being calculated. Proposed algorithm however may will give further lower run times for certain values in $F_p^*$
3. Proposed method can be  modified with same run time advantages for Elliptic curves too.

\section{Conclusion and Acknowledgements}
Discrete Logarithm problem  is classical problem , cryptographic systems based on hardness of this  problem still remain very popular. However special cases  where p is smooth prime proposed method provides practically faster algorithm than Pohlig-Hellman . However reducing the upper bound of the worst case still remain a challenge. \\
 
I would like to thank Prof G.P.Biswas of IIT Dhanbad and Prof R.Munshi of ISI Kolkata, for their constant support and guidance. I owe special thanks to Gaurav sinha, IRS , and Syed Waquar Raza,IPS, for reading  the proof of the paper and giving valuable suggestions.

\includegraphics[scale=0.75,page=1]{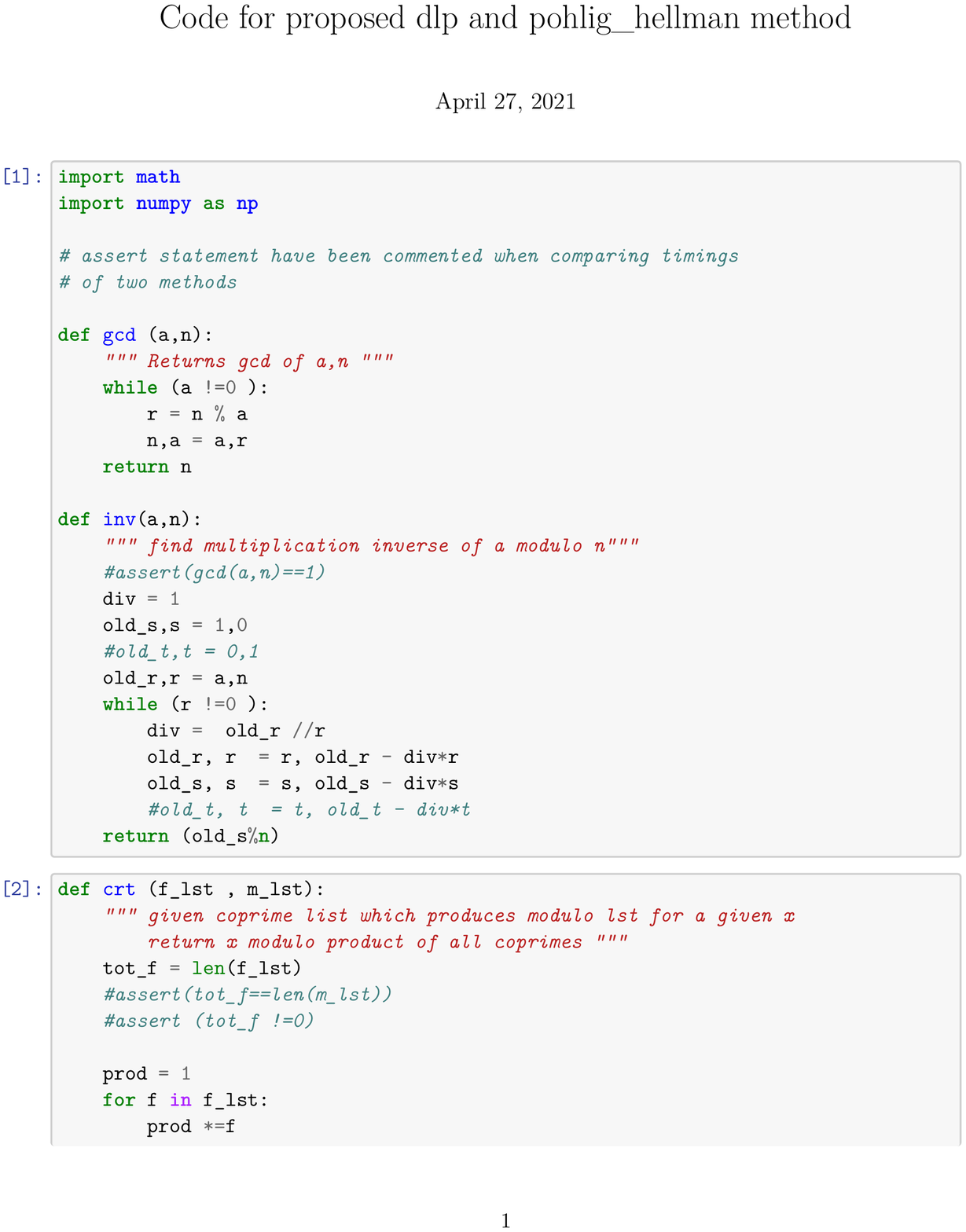}
\includegraphics[scale=0.75,page=2]{code.pdf}
\includegraphics[scale=0.75,page=3]{code.pdf}
\includegraphics[scale=0.75,page=4]{code.pdf}
\includegraphics[scale=0.75,page=5]{code.pdf}
\end{document}